\documentclass[11pt]{article}
\usepackage{amssymb, amsmath, amsfonts, latexsym, hyperref}
\usepackage{enumerate,pspicture,epsfig, amsthm}
\usepackage{graphicx,mathtools}
\usepackage[font=small,labelfont=bf]{caption}

\usepackage[usenames,dvipsnames]{pstricks}
\usepackage{pst-grad, pst-plot}
\usepackage{hyperref, float, authblk}
\usepackage[english]{babel}
\usepackage{caption}
\everymath{\displaystyle}
\usepackage[margin=1in]{geometry}
\usepackage{multirow,booktabs}
\usepackage[mathscr]{euscript}
\floatstyle{boxed}
\title{Betweenness Centrality of Cartesian Product of Graphs}
\author[1]{Sunil Kumar R}
\author[2]{Kannan Balakrishnan}
\affil[1,2]{Department of Computer Applications, Cochin University of Science and Technology}
\affil[1]{sunilstands@gmail.com}
\affil[2]{mullayilkannan@gmail.com}
          \date{ }           
             
\newtheorem{thm}{Theorem}[section]
\newtheorem{lem}{Lemma}[section]

\newtheorem{cor}{Corollary}[section]
\newtheorem{note}{Note}[section]
\newtheorem{prop}{Proposition}[section]
\newtheorem{defn}{Definition}[section]

\begin{document}
\maketitle

\begin{abstract} 
Betweenness centrality is a widely-used measure in the analysis of large complex networks. It measures the potential or power of a vertex to control the communication over the network under the assumption that information primarily flows over the shortest paths between them.  In this paper we prove several results on betweenness centrality of Cartesian product of graphs.
\vspace{.2cm}
\end{abstract} 
 \section{Introduction}
 Several Centrality measures have so far been studied and their importance is increasing day by day. \textit{Betweenness centrality} has a vital role in the  analysis of networks.\cite{de2015dynamic,martin2010centrality,rubinov2010complex,salathe2010high}  
It has many applications in a variety of domains such as biological networks \cite{jeong2001lethality,pinney2005decomposition,del2005topology,koschutzki2008centrality,estrada2006virtual}, study of sexual networks and AIDS\cite{liljeros2001web}, identifying key actors in terrorist networks\cite{krebs2002mapping}, transportation networks\cite{guimera2005worldwide}, supply chain management\cite{cisic2000research}, bio-informatics-protein interaction networks\cite{joy2005high,chen2009disease}, food webs \cite{jordan2009keystone}etc. Betweennes centrality \cite{borgatti2006graph,brandes2008variants} indicates the betweenness of a vertex (or an edge) in a network and it measures the extent to which a vertex (or an edge) lies on the shortest paths between pairs of other vertices.  It is quite difficult to find out the betweenness centrality of a vertex in a large graph. The  computation of this index based on direct application of definition becomes impractical as the number of nodes $n$ increases and has complexity in the order of $\mathcal{O}(n^3)$. The fastest exact algorithm due to Brandes\cite{brandes2001faster} requires $\mathcal{O}(n+m)$ space and $\mathcal{O}(nm)$ time where $n$ is the number of nodes and $m$ the number of edges in the graph.  Exact computations of betweenness centrality can take a lot of time even for Brandes algorithm. But a large network can be thought of as it is made by joining smaller networks together. There are several graph operations which results in a larger graph $G$ and many of the properties of larger graphs can be derived from its constituent graphs. Graph operations are used for constructing new classes of graphs. Cartesian product is an important graph operation.\\

It is assumed that the graphs taken here are simple undirected connected graphs. 
Graph-reference may be given for making the context clear.
\section{Background}
  The concept of betweenness centrality of a vertex was first introduced by Bavelas in 1948 \cite{bavelasa}. The 
importance of the concept of vertex centrality is that how a vertex  acts as a bridge among all the pairs of vertices in joining them by shortest paths. It gives the potential of a vertex for control of information flow in the network\cite{brandes2005centrality,borgatti2005centrality}. The \textit{order} of a graph $G$ is the number of vertices in $G$; it is denoted by $∣|G|∣$. The same notation is used for the number of elements (cardinality) of a set. Thus, $∣|G|∣=∣|V(G)|∣$. The distance between
two vertices $u, v \in V (G)$, denoted by $d_ G (u, v)$, is the length of the shortest path in $G$ between $u$ and $v$.  A shortest path joining vertices $u$ and $v$  is called a \textit{geodesic} between $u$ and $v$. A graph $G$ is a \textit{geodetic graph}\cite{ore1962theory} if every pair of vertices of $G$ is connected by a unique shortest path.
 The diameter, $d iam(G)$, of a graph $G$ is given by $\max\{d(u, v) | u, v \in V (G)\}$. Two vertices $u$ and $v$ of $G$ with $d(u,v)=diam(G)$ are \textit{diametrical vertices}\cite{mulder1980interval}.  The interval $I_G(u, v)$ consists of all vertices on geodesics joining $u$ and $v$ in $G$.

A graph $G$ is \textit{vertex-transitive} if every vertex in $G$ can be mapped to any other vertex by some automorphism. Similarly a graph is \textit{edge-transitive} if its automorphism group acts transitively on the set of edges. 

A definition to betweenness centrality of a vertex in a graph $G$, given by Freeman \cite{freeman1977set} is as follows
\begin{defn}[Betweenness Centrality]%\cite{freeman1977set}
If $x\in V(G)$, the betweenness centrality $B(x)$ for $x$ is defined as 
$$B(x) = \sum_{u\neq v\neq x}{\delta(u,v|x)}$$
provided $\delta(u,v|x)=\frac{\sigma(u,v|x)}{\sigma (u,v)}$ where $ \sigma (u,v)$ is the number of shortest $u$-$v$ paths  and $\sigma(u,v|x)$ is the number of shortest $u$-$v$ paths containing    $x$. The ratio $\delta(u,v|x)$ is called the pair-dependancy of the pair of vertices $\{u,v\}$ on $x$.\end{defn}
Observe that $x \in V(G)$ lies on the shortest path between two vertices $u, v \in V(G)$ , iff $d(u, v) = d(u, x)+d(x, v)$. The number of shortest $u$-$v$ paths passing through $x$ is given by
\begin{equation}\label{eq1}
\sigma(u,v|x) = \sigma(u,x)\times \sigma(x,v)
\end{equation} 
\begin{defn}[Cartesian Product, \cite{sabidussi1957graphs}]
%\subsection{Cartesian Product }
The \textit{Cartesian product} of two graphs $G$ and $H$, denoted by $G \square H$, is a graph with vertex set $V (G)\times V (H)$, where two vertices $(g, h)$ and $(g' , h' )$ are adjacent if $g = g'$ and $hh' \in E(H)$, or $gg' \in E(G)$ and $h = h'$.     The graphs $G$ and $H$ are called \textit{factors} of the product $G \square H$.\end{defn} For any $h \in V (H)$, the subgraph of $G\square H$ induced by $V (G) \times \{h\}$  is called as \textit{$G$-fiber} or \textit{$G$-layer}, denoted by $G^h$ . Similarly, we can define \textit{$H$-fiber} or \textit{$H$-layer}. They are isomorphic to $G$ and $H$, respectively. $G\square H$ contains $|H|$ copies of $G$ and $|G|$ copies of $H$.  Projections are the maps from a product graph to its factors. They are weak homomorphisms in the sense that they respect adjacency. The two projections on $G\square H$ namely $p_G:G \square H \rightarrow G$ and $p_H:G \square H \rightarrow H$ defined by $p_G(g,h)=g$ and $p_H(g,h)=h$ refer to the corresponding $G$-, $H$- coordinates. Thus an edge in $G \square H$ is mapped into a single  vertex by one of the projections $p_G$ or $p_H$ and into an edge by the other. If $G$ and $H$ are connected, then
$G\square H$ is also connected. Assuming isomorphic graphs are equal, Cartesian product is commutative as well as associative. \\

For a graph $G$ and $v \in V (G)$, the degree of a vertex $v$ is denoted by $d_G (v)$, or simply
$d(v)$. Furthermore, we denote by $\delta(G)$ the
minimum degree of a graph $G$. The minimum degree is additive under Cartesian products,
i.e. $\delta(G\square H) = \delta(G) + \delta(H)$. Recall that the symbol $N_G (v)$ denotes the set of neighbours of a vertex $v$ in a graph $G$. Thus $d_G(v)=|N_G (v)|$.\\
\begin{defn}[Cartesian product of several graphs,\cite{imrich2008topics}]
The \textit{Cartesian product} $G=G_1\square G_2\square \ldots \square G_k$ of the graphs $G_1,G_2,\ldots,G_k$ is defined on the $k$-tuples $(v_1,v_2,\ldots,v_k)$, where $v_i\in G_i,\; 1\leq i\leq k$ in such a way that two $k$-tuples $(u_1,u_2,\ldots,u_k)$ and $(v_1,v_2,\ldots,v_k)$ are adjacent if there exists an index  $l$ such that $[u_l,v_l]\in E(G_l)$ and $u_i=v_i$ for $i\neq l$. The $k$-tuples $(v_1,v_2,\ldots,v_k)$ are called \textit{coordinate vectors}, and the $v_i$ are the \textit{coordinates}.
\end{defn}
The Cartesian product $G=G_1\square G_2\square \ldots \square G_k$ of $k$-factors is briefly denoted as $G=\square_{i=1}^kG_i$. The $n^{th}$ Cartesian product of a graph $G$ is  donoted as $G^n=\square_{i=1}^nG$. It is to be noted that the product $G$ is connected if and only if each of its factor $G_i$ is connected and the diameter of the product is given by, $diam(\square_{i=1}^kG_i)=\sum_{i=1}^kdiam(G_i)$. \\
 
The following proposition shows that the distance between two vertices in the product graph is the sum of the distance between their projections in the factor graphs.
\begin{lem}\label{a}\cite{hammack2011handbook}
If $(g, h)$ and $(g' , h' )$ are vertices of a Cartesian product $G \square H$, then\\
\begin{equation}
d_{G \square H} \big[(g, h), (g' , h' )\big] = d_G (g, g' ) + d_H (h, h' )
\end{equation}
\end{lem}

This can be generalized to the following lemma.
\begin{lem}[Distance lemma]\label{a1}\cite{hammack2011handbook}
Let $G$ be the Cartesian product $G=\square_{i=1}^k G_i$ of connected graphs, and let $g=(g_1,\ldots,g_k)$ and $g'=(g_1',\ldots,g_k')$ be vertices of $G$. Then
$$d_G(g,g')=\sum_{i=1}^k d_{G_i}(g_i,g'_i)$$
\end{lem}
 Lemma \ref{a} implies that $d_{G \square H} \big[(g, h), (g',h)\big] = d_{G^h} (g, h), (g',h)$. In other words, $d_{G \square H}$ restricted to $G^h$ is $d_{G^h}$. It means that every shortest path in a $G$-fiber ia also a shortest path in  ${G \square H}$. Subgraphs with this property are called \textit{isometric}. That is, a subgraph $U$ of a graph $G$ is isometric in $G$ if $d_U(u,v)=d_G(u,v)$ for all $u,v\in G$. It can be easily seen that $G$-fibers ($H$-fibers) are isometric subgraphs of ${G \square H}$. Every shortest ${G \square H}$-path between two vertices of one and the same fiber $G^h$ or $H^g$ is already in that fiber. Such subgraphs are called \textit{convex}. A subgraph $U$ of a graph $G$ is convex in $G$ if every shortest $G$-path between vertices of $U$ is already in $U$.
 \begin{lem}\cite{imrich2008topics}
 Let $G$ and $H$ be connected graphs. Then all $G$-fibers and $H$-fibers are convex subgraphs of $G\square H$.
 \end{lem}
 
For a connected graph $G$ and $u,v\in G$, the \textit{interval} $I_G(u,v)$ between $u$ and $v$ is defined as the set of vertices that lie on shortest $u$-$v$ paths; that is,\\
$$I_G(u,v)=\{w \in G:\;d(u,v)=d(u,w)+d(w,v)\}$$
\begin{prop}\cite{imrich2008topics}
Let $v_1=(g_1, h_1)$ and $v_2=(g_2 , h_2)$ be two vertices of  $G \square H$, then the vertex  $v_3=(g_3,h_3)$ lies in $I_{G \square H}(v_1,v_2)$ if and only if $g_3 \in I_G(g_1,g_2)$  and $h_3 \in I_H(h_1,h_2)$.
\end{prop}
It can be generalized as follows.
\begin{prop}
Let $G=\square_{i=1}^k G_{n_i}$. Let $g(g_1,g_2,\ldots, g_k)$, $g'(g'_1,g'_2,\ldots,g'_k)$ and $ g''(g''_1,g''_2,\ldots,g_k'')$ are any three vertices in $G$. Then $g''$ lies in the shortest path of $g$ and $g'$ if and only if $g''_i \in I(g_i,g'_i)$ $\forall i$.
\end{prop}
\subsection{The betweenness centrality of vertices in Cartesian product of two graphs }

The following proposition shows how the number of geodesics between two vertices $u$ and $v$ in a product graph $G\square H$ is related to the the number of geodesics between their projections in the factor graphs.
\begin{prop}\label{j}
If $u=(g, h)$ and $v=(g', h')$ are vertices of  $G \square H$, then the number of shortest paths, $\sigma_{G \square H}$, between them in $G \square H$ is given by
\begin{equation}\label{e1}
\sigma_{G \square H}\big[(g,h),(g', h' )\big]=\sigma_G(g, g')\times \sigma_H(h, h' )\times \binom{d_G(g, g')+d_H(h,h')}{d_G (g,g')}
\end{equation}
\end{prop}
\begin{proof}
Consider the vertices $u=(g, h)$ and $v=(g', h')$ in $G\square H$. Let $d=d(u,v)$ denote the distance between $u$ and $v$ in $G \square H $. Suppose there exists unique shortest paths between $g,g'$ and $h,h'$. Every shortest path from $u$ to $v$ is a sequence of $d$ edges and the image of each edge is an edge lying between $g$ and $g'$ or $h$ and $h'$ under the projections $p_G$ and $p_H$. Let a sequence of $d$ edges in the $u$-$v$ path in the product makes a sequence of $d_G$ edges in $G$ and a sequence of $d_H$ edges in $H$ so that $d=d_G+d_H$. Since $d_G$ and $d_H$ are the same for any $u$-$v$ path, the number of shortest paths  between $u$ and $v$ is the number of ways of selecting $d_G$ edges from $d$ edges, which is $\binom{d}{d_G}$. If there exists $\sigma _G$ shortest paths between $g$ and $g'$ in $G$ and $\sigma _H$ shortest paths between $h$ and $h'$ in $H$, then corresponding to each pair there exists $\binom{d}{d_G}$  shortest paths between $u$ and $v$ in $G \square H$. \\Therefore,  $\sigma_{G \square H}\big[(g,h),(g', h' )\big]=\sigma_G(g, g')\times \sigma_H(h, h' )\times \binom{d_G(g, g')+d_H(h,h')}{d_G (g,g')}$\\
For brevity, we may write  $\sigma= \sigma_G \times \sigma_H \times \binom{d}{d_G}$
\end{proof}
\begin{cor}
If $G$ and $H$  are geodetic graphs, then the number of shortest paths between $(g,h)$ and $(g',h')$ in $G\square H$ is given by
$$\sigma_{G \square H}\big[(g,h),(g', h' )\big]=\binom{d_G(g, g')+d_H(h,h')}{d_G (g,g')}$$ 
\end {cor}
By the associativity of $\square$, equation \ref{e1} can be generalized as 
\begin{prop}\label{d}
Let $G=\square_{i=1}^k G_{n_i}$. If $u=(u_1,u_2,\ldots,u_k)$, $v=(v_1,v_2,\ldots,v_k)$ are two vertices in $G$ such that
 $\sigma_{G_i}(u,v)=\sigma_i$, $d_{G_i}(u_i,v_i)=d_i$ and $d=\sum d_i$, then\\
$$\sigma_{G}(u,v)=\sigma_{1} \sigma_{2} \ldots \sigma_{n} \binom{d}{d_1} \binom{d-d_1}{d_2} \binom{d-d_1-d_2}{d_3} \ldots 1 $$
\end{prop}
\begin{prop}\label{c}
Let $v_1$,$v_2$ and $v_3$ are any three vertices in $G \square H$. Then\\
\begin{equation}\label{eq4}
\sigma_{G \square H}(v_1,v_2|v_3)=\sigma_{G \square H}(v_1,v_3)\times\sigma_{G \square H}(v_3,v_2)
\end{equation}
\end{prop}
\begin{thm}
If $u=(g,h)$, $v=(g',h')$ be any  distinct vertices in $G\square H$ then the betweenness centrality of  $x=(g_0,h_0)$ in $G\square H$ is given by\\ 
$$B_{G\square H}(x) =\sum_{u\neq v\neq x}{\delta_{G\square H}(u,v|x)}$$ where
$$\delta_{G\square H}(u,v|x)=\delta_{G}(g,g'|g_0)\times \delta_{H}(h,h'|h_0))\times \frac{d_1\times d_2}{d}$$ where
$$d_1=\binom{d_{G}(g,g_0)+d_{H}(h,h_0)}{d_G(g,g_0)},\;d_2={\binom{d_{G}(g_0,g')+d_{H}(h_0,h')}{d_G(g_0,g')}}, \;d={\binom{d_{G}(g,g')+d_{H}(h,h')}{d_G(g,g')}}$$\\
\end{thm}
\begin{proof}
The result follows from the definition of betweenness centrality and from equations \ref{eq1}-\ref{eq4}\\ Hence $$\delta_{G\square H}(u,v|x) = \frac{\sigma_{G\square H}(u,v|x)}{\sigma_{G\square H} (u,v)}\\=\frac{\sigma_{G}(g,g'|g_0)}{\sigma_{G} (g,g')}\times\frac{\sigma_{H}(h,h'|h_0)}{\sigma_{H} (h,h')}\times \frac{d_1\times d_2}{d}$$\\
\end{proof}
\section{Wiener index of a graph}
The Wiener index \cite{wiener1947structural} of a graph $G$, denoted by $W(G)$ is the sum of the distances between all (unordered)  pairs of vertices of $G$. That is, \\
$$W(G)=\sum_{i<j}d(v_i,v_j)$$
or,
$$W(G)=\frac{1}{2}\sum_{u,v \in V(G)}
d(u,v)$$
Wiener index is also named \textit{total status} or \textit{total distance} of a graph. The Wiener index of Cartesian product \cite{graovac1991wiener,yeh1994sum,dehmer2014quantitative} of two graphs $G$ and $H$ is given by\\
\begin{equation}
W(G\square H)=|G|^2W(H)+|H|^2W(G)
\end{equation}
It can be extended to\\
\begin{equation}\label{eq4}
W(\square_{i=1}^nG_i)=\sum_{i=1}^n\Big(W(G_i)\prod_{j\neq i}|G_j|^2\Big)
\end{equation}
The betweenness centrality of $G$ is given by\\
\begin{equation}\label{eq2}
\sum _{v\in V(G)}B(v)=W(G)-\binom{|G|}{2}
\end{equation}
The average distance of a graph $G$, denoted by $\mu(G)$ is given by\\
$$\mu(G)=\frac{\hbox{Total distance}}{\hbox{No. of distinct pairs}}=\frac{W(G)}{\binom{|G|}{2}}$$
\subsubsection{Grid graphs}
\textit{Grid  graphs} are the cartesian product of path graphs. $P_m\square P_n$ represents a rectangular grid $R$. If $u=(g,h)$ and $v=(g',h')$ are any two vertices of $R$, then $d(u,v)=|g-g'|+|h-h'|$ and $\sigma (u,v)=\binom{d}{d_1}$ where $d_1=|g-g'|$ or $|h-h'|$.
 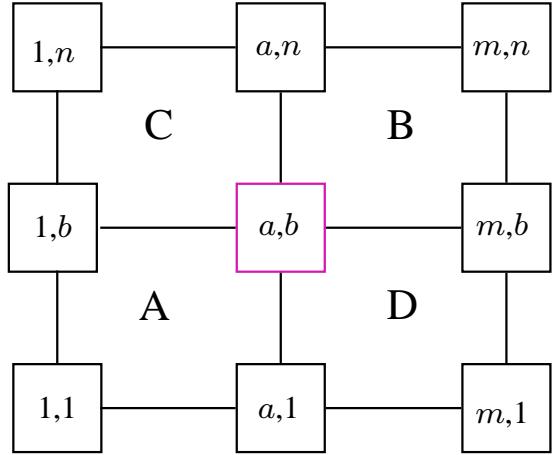
\begin{figure}[H]
\centering
\scalebox{1.5}
{
\begin{pspicture}(0,-2.0)(6.0474997,2.0)
\definecolor{color457}{rgb}{0.8274509803921568,0.09411764705882353,0.6901960784313725}
\usefont{T1}{ptm}{m}{n}
\rput(2.0414062,0.91){C}
\usefont{T1}{ptm}{m}{n}
\rput(4.1870313,0.91){B}
\usefont{T1}{ptm}{m}{n}
\rput(2.0023437,-0.69){A}
\usefont{T1}{ptm}{m}{n}
\rput(4.201406,-0.69){D}
\usefont{T1}{ptm}{m}{n}
\rput(1.1,1.55){\centering \scriptsize 1,$n$}
\usefont{T1}{ptm}{m}{n}
\rput(3.1,1.57){\centering \scriptsize $a$,$n$}
\usefont{T1}{ptm}{m}{n}
\rput(5.078906,1.57){\centering \scriptsize $m$,$n$}
\usefont{T1}{ptm}{m}{n}
\rput(1.10,-0.01){\centering \scriptsize 1,$b$}
\usefont{T1}{ptm}{m}{n}
\rput(3.1,0.01){{\centering \scriptsize $a$,$b$}}
\usefont{T1}{ptm}{m}{n}
\rput(5.08375,0.01){\centering \scriptsize $m$,$b$}
\usefont{T1}{ptm}{m}{n}
\rput(1.1415623,-1.61){\centering \scriptsize 1,1}
\usefont{T1}{ptm}{m}{n}
\rput(3.103125,-1.61){\centering \scriptsize $a$,1}
\usefont{T1}{ptm}{m}{n}
\rput(5.082344,-1.66){\centering \scriptsize $m$,1}
\psframe[linewidth=0.02,dimen=outer](3.520781,-1.2)(2.720781,-2.0)
\psframe[linewidth=0.02,linecolor=color457,dimen=outer](3.520781,0.4)(2.720781,-0.4)
\psframe[linewidth=0.02,dimen=outer](3.520781,2.0)(2.720781,1.2)
\psframe[linewidth=0.02,dimen=outer](5.520781,0.4)(4.7207813,-0.4)
\psframe[linewidth=0.02,dimen=outer](5.520781,2.0)(4.7207813,1.2)
\psframe[linewidth=0.02,dimen=outer](5.520781,-1.2)(4.7207813,-2.0)
\psline[linewidth=0.02cm](1.5207812,1.6)(2.720781,1.6)
\psline[linewidth=0.02cm](3.520781,1.6)(4.7207813,1.6)
\psline[linewidth=0.02cm](1.5207812,-1.6)(2.720781,-1.6)
\psline[linewidth=0.02cm](3.520781,-1.6)(4.7207813,-1.6)
\psline[linewidth=0.02cm](3.1207812,-0.4)(3.1207812,-1.2)
\psline[linewidth=0.02cm](3.1207812,1.2)(3.1207812,0.4)
\psline[linewidth=0.02cm](3.520781,0.0)(4.7207813,0.0)
\psline[linewidth=0.02cm](1.5207812,0.0)(2.720781,0.0)
\psline[linewidth=0.02cm](5.120781,1.2)(5.120781,0.4)
\psline[linewidth=0.02cm](5.120781,-0.4)(5.120781,-1.2)
\psframe[linewidth=0.02,dimen=outer](1.5407811,2.0)(0.7407814,1.2)
\psframe[linewidth=0.02,dimen=outer](1.500781,0.4)(0.7007814,-0.4)
\psframe[linewidth=0.02,dimen=outer](1.5407811,-1.2)(0.7407814,-2.0)
\psline[linewidth=0.02cm](1.14,1.22)(1.14,0.38)
\psline[linewidth=0.02cm](1.14,-0.38)(1.14,-1.2)
\end{pspicture}

} 
\caption{Rectangular grid $P_m\square P_n$}\label{fig5}
\end{figure}
Consider the rectangular grid $P_m\square P_n$.  Let $x$ denotes the vertex $(a,b)\in P_m\square P_n$ where $1\leq a\leq m,\;1\leq b\leq n$.  Consider the paths $a\square P_n$ and $b\square P_m$ passing through $(a,b)$. They divide the rectangular grid into 4 parts namely $A,B,C,D$ sharing their common sides. Figure\ref {fig5}. Now the pairs of vertices in the diagonal regions $A,B$ and $C,D$ contribute to the betweenness centrality of $(a,b)$. Hence \\
$$B\big[(a,b)\big]=  \sum_{u,v}\frac{\sigma(u,x)\sigma(x,v)}{\sigma(u,v)}-\Big[(a-1)\times(m-a)+(b-1)\times(n-b)\Big]$$
 where $u$ and $v$ belongs the diagonal quadrants $A,B$ and then $C,D$.
 \begin{note}
 The number of geodesics of length $kn$ in the grid $G=\square_{i=1}^{k}P_{n+1}$ from $(0,0,\ldots,k\; times)$ to $(n,n,\ldots k\; times)$ can be obtained from the $k$-ary de Bruijn sequence $s(k,n)$ where $s(k,n)=\frac{(kn)!}{(n!)^k}$ ($OEIS$ $A000984$)
\end{note}
\subsection{Hamming graphs}
\textit{Hamming graphs} are Cartesian products of complete graphs. If $G$ is a Hamming graph, then  
$G=K_{n_1}\square K_{n_2}\square\dots  \square K_{n_r}$ for some $r\geq 1$ and $n_i\geq 2$. The vertices of $G$ can be labeled with vector $(a_1,a_2,\ldots a_r)$ where $a_i \in \{0,1,\dots,n_i-1\}$. Two vertices of $G$ are adjacent if the corresponding tuples differ in precisely one coordinate.
The distance (named  \textit{Hamming distance}) between two vertices $u$ and $v$ denoted by $d(u,v)$ is the number of positions in which the two vectors differ. 

\textit{Hypercubes} are Cartesian product of complete graphs $K_2$. An $r$-dimensional hypercube (or $r$-cube) denoted by $Q_r$ is given by, $Q_r=\square_{i=1}^{r}K_2$. It can also be defined recursively, $Q_r=K_2\square Q_{r-1}$. Hypercubes are important classes of graphs having many interesting structural properties. The number of geodesics between $u,v\in Q_r$ is given by $\sigma(u,v)=d(u,v)!$. For a connected graph $G$, the condition ``$I(u,v)$ induces a $d(u,v)$-dimensional hypercube for any two vertices $u$ and $v$ of $G$ implies that $G$ " is a Hamming graph \cite{mulder1980interval}.
\begin{lem}\cite{peterin2002characterizing}
A graph $G$ is a nontrivial subgraph of the Cartesisn product of graphs if and only if $G$ is a nontrivial subgraph of the Cartesian product of two complete graphs.
\end{lem}
\begin{prop}
Let $G$ be the Hamming graph $G=K_{n_1}\square K_{n_2}\square\dots  \square K_{n_r}$. Then the betweenness centrality of $v\in G$ is given by\\
\begin{equation}
B(v)=\frac{1}{2}\prod_{i=1}^{r}n_i\Big[r-1-\sum_{i=1}^r\frac{1}{n_i}\Big]+\frac{1}{2}
\end{equation}
\end{prop}\label{z}

\begin{proof}
Let $H=\square_{i=1}^r K_{n_i}$. Since $H$ is vertex transitive,
from equations \ref{eq4} and \ref{eq2},
\begin{align*}
W(H)&=\sum_{i=1}^r W(K_{n_i})\prod_{j=1,j\neq i}^r |K_{n_j}|^2 
 =\sum_{i=1}^r\binom{n_i}{2}\prod_{j=1,j\neq i}^r n_j^2\\
 &=\frac{1}{2}\Big(\prod_{i=1}^r n_i\Big)^2\Big[r-\sum_{i=1}^r\frac{1}{n_i}\Big]\\
 B_H(v)&=\frac{W(H)-\binom{|H|}{2}}{|H|}\\
 &=\frac{1}{2}\prod_{i=1}^{r}n_i\Big[r-1-\sum_{i=1}^r\frac{1}{n_i}\Big]+\frac{1}{2}
\end{align*}
\end{proof}
\begin{cor}
If $K_p$, $K_q$ and $K_r$ are complete graphs, then\\\\ for $v \in K_p \square K_q$
$$B(v)=\frac{(p-1)(q-1)}{2}$$ 
for $v\in K_p \square K_q \square K_r$\\
$$B(v)=\frac{1}{2}\Big[2pqr-(pq+pr+qr)+1\Big]$$ 
\end{cor}
\begin{cor}
If $v \in \square_{i=1}^r K_{n}$, then\\
$$B(v)=\frac{1}{2}\Big[(r-1) n^r -r n^{r-1}+1\Big]$$
when $n=2$, $v \in \mathrm{Q_r}$, the $r$-cube, then\\
$$B(v)=(r-2)2^{r-2}+\frac{1}{2}$$
\end{cor}
\begin{figure}[H]
\centering
\scalebox{1} % Change this value to rescale the drawing.
{
\begin{pspicture}(0,-2.94)(5.14,2.94)
\psframe[linewidth=0.04,dimen=outer](3.52,1.4)(1.62,-0.5)
\psframe[linewidth=0.04,dimen=outer](3.86,1.74)(1.22,-0.9)
\psframe[linewidth=0.04,dimen=outer](4.24,2.08)(0.84,-1.32)
\psframe[linewidth=0.04,dimen=outer](4.68,2.48)(0.44,-1.76)
\psframe[linewidth=0.04,dimen=outer](5.06,2.86)(0.04,-2.16)
\psline[linewidth=0.04cm](3.52,1.36)(5.06,2.88)
\psline[linewidth=0.04cm](1.66,-0.5)(0.1,-2.1)
\psline[linewidth=0.04cm](1.66,1.36)(0.08,2.82)
\psline[linewidth=0.04cm](3.46,-0.46)(5.06,-2.16)
\psbezier[linewidth=0.04](3.52,1.32)(3.52,0.52)(5.02,2.04)(5.02,2.84)
\psbezier[linewidth=0.04](3.5,-0.48)(3.5,-1.28)(5.02,-2.92)(5.02,-2.12)
\psbezier[linewidth=0.04](1.66,-0.46)(1.66,-1.26)(0.06,-2.92)(0.06,-2.12)
\psbezier[linewidth=0.04](1.68,1.32)(1.68,0.52)(0.08,2.02)(0.08,2.82)
\psdots[dotsize=0.16](1.28,1.68)
\psdots[dotsize=0.16](1.7,1.32)
\psdots[dotsize=0.16](0.92,2.02)
\psdots[dotsize=0.16](0.5,2.42)
\psdots[dotsize=0.16](0.08,2.78)
\psdots[dotsize=0.16](0.08,-2.12)
\psdots[dotsize=0.16](0.5,-1.74)
\psdots[dotsize=0.16](0.9,-1.32)
\psdots[dotsize=0.16](1.28,-0.86)
\psdots[dotsize=0.16](1.66,-0.5)
\psdots[dotsize=0.16](3.52,-0.5)
\psdots[dotsize=0.16](3.84,-0.86)
\psdots[dotsize=0.16](4.26,-1.3)
\psdots[dotsize=0.16](4.68,-1.72)
\psdots[dotsize=0.16](5.04,-2.16)
\psdots[dotsize=0.16](5.02,2.84)
\psdots[dotsize=0.16](4.68,2.46)
\psdots[dotsize=0.16](4.26,2.08)
\psdots[dotsize=0.16](3.88,1.7)
\psdots[dotsize=0.16](3.54,1.32)
\end{pspicture} 
}
\caption{$C_4\square C_5$}
\end{figure}
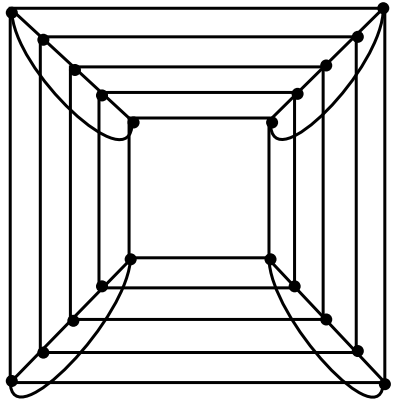
\subsubsection{Product of cycles}
 For the cycle $C_n$, $W(C_n)=\frac{1}{8}n^3$ when $n$ is even, and $W(C_n)=\frac{1}{8}(n^3-n)$ when $n$ is odd. Product of cycles $\square_{i=1}^r C_{n_i}$ is vertex transitive. Therefore from equations \ref{eq4} and\ref{eq2} we get\\\\
for $n_i\in 2\mathbb{Z}^+$ $$W(\square_{i=1
}^n C_{n_i})==\frac{1}{8}\Big(\prod_{i=1}^r n_i\Big)^2\;\sum_{i=1}^r n_i$$ \\ and for $n_i\in 2\mathbb{Z}^++1$ $$W(\square_{i=1
}^n C_{n_i})==\frac{1}{8}\Big(\prod_{i=1}^r n_i\Big)^2\;\sum_{i=1}^r \Big(n_i-\frac{1}{n_i}\Big)$$ 
\begin{prop}
If $G$ is the Cartesian product of $r$ even cycles. \\
ie. $G=\square_{i=1}^r C_{n_i}$ , $n_i\in 2\mathbb{Z}^+$, then for $v \in G$\\
$$B(v)=\frac{1}{8}\Big[\prod_{i=1}^rn_i\;\sum_{i=1}^rn_i-4\Big(\prod_{i=1}^rn_i-1\Big)\Big]$$
if $n_i=2k_i$,
$$B(v)=2^{r-2} \prod_{i=1}^r k_i \; \Big[\sum_{i=1}^r k_i-2 \Big]+\frac{1}{2}$$
\end{prop}
\begin{prop}
If $G$ is the Cartesian product of $r$ odd cycles. \\
ie. $G=\square_{i=1}^r C_{n_i}$ , $n_i\in 2\mathbb{Z}^++1$, then for $v \in G$\\
$$B(v)=\frac{1}{8}\Big[\prod_{i=1}^rn_i\;\sum_{i=1}^r\Big(n_i-\frac{1}{n_i}\Big)-4\Big(\prod_{i=1}^rn_i-1\Big)\Big]$$
\end{prop}
\begin{cor}
Consider two cycles $C_m$ and $C_n$. Let $v\in C_{m}\square C_{n}$ then\\ 
$$B(v)= \begin{cases}\frac{1}{8}(mn-1)(m+n-4)\;,\; when\;m\;and\;n\;are\; odd\\
       \frac{1}{8}(mn^2+m(m-4)n+4)\qedhere\;,\; when\;m\;and\;n\;are\; even\\
       \frac{1}{8}[mn^2+(m^2-4m-1)n+4]\; ,\;when\; m\;odd\;and\;n\; is\; even.\qedhere       
       \end{cases}$$
In another form,\\

$$B(v)= \begin{cases}k_1k_2(k_1+k_2)+\binom{k_1}{2}+\binom{k_2}{2} \;\;\hbox{when  }        			m=2k_1+1,n=2k_2+1 \\    
k_1k_2(k_1+k_2-2)+\frac{1}{2} \;\;\hbox{when  }m=2k_1,n=2k_2\\       
k_1k_2(k_1+k_2-1)+\frac{1}{2} (k_2-1)^2 \;\;\hbox{when  }m=2k_1+1,n=2k_2 \\  
\end{cases}        
$$
\end{cor}
 \section{Conclusion}
A composite graph can be constructed by applying different graph oparations on smaller graphs and hence many of the structural properties of the composite graphs can be studied from its constituent smaller subgraphs. Here we tried to find the betweenness centrality of Cartesian product of graphs. This can be extended to other products also.
\bibliographystyle{unsrt}
\bibliography{sunilbiblio}
\end{document}